\documentclass[titlepage,draft,12pt]{article} 
\usepackage{amsfonts,latexsym} 
\usepackage[a4paper]{geometry}

\date{}

\newcommand{\ep}{\varepsilon}
\newcommand{\qed}{{\penalty 10000\mbox{$\quad\Box$}}}
\newcommand{\re}{\mathbb{R}}

\newcommand{\n}{\mathbb{N}}
\newcommand{\cep}{c_{\ep}}
\newcommand{\Cep}{C_{\ep}}
\newcommand{\gep}{g_{\ep}}
\newcommand{\rep}{r_{\ep}}
\newcommand{\roep}{\rho_{\ep}}
\newcommand{\uepn}{u_{\ep,\nu}}
\newcommand{\uepm}{u_{\ep,\mu}}
\newcommand{\psiep}{\psi_{\ep}}
\newcommand{\uep}{u_{\ep}}
\newcommand{\yep}{y_{\ep}}
\newcommand{\tetep}{\theta_{\ep}}
\newcommand{\auqg}[1]{|A^{1/2}#1|^{2\gamma}}
\newcommand{\dat}{D(A^{3/2})}
\newcommand{\da}{D(A)}
\newcommand{\dau}{D(A^{1/2})}

\newtheorem{thm}{Theorem}[section]
\newtheorem{thmbibl}{Theorem}
\newtheorem{propbibl}[thmbibl]{Proposition}
\newtheorem{rmk}[thm]{Remark}
\newtheorem{prop}[thm]{Proposition}
\newtheorem{defn}[thm]{Definition}

\newtheorem{lemma}[thm]{Lemma}

\title{Optimal decay-error estimates for the hyperbolic-parabolic
singular perturbation of a degenerate nonlinear equation}

\author{Marina Ghisi\vspace{1ex}\\ 
{\normalsize Universit\`a degli Studi di Pisa} \\
{\normalsize Dipartimento di Matematica ``Leonida Tonelli''}\\ 
{\normalsize PISA (Italy)}\\
{\normalsize e-mail: \texttt{ghisi@dm.unipi.it}}
\and
Massimo Gobbino\vspace{1ex}\\ 
{\normalsize Universit\`a degli Studi di Pisa} \\
{\normalsize Dipartimento di Matematica Applicata ``Ulisse Dini''}\\ 
{\normalsize PISA (Italy)}\\  
{\normalsize e-mail: \texttt{m.gobbino@dma.unipi.it}}}

\begin{document}
\maketitle
\begin{abstract}
	We consider a \emph{degenerate} hyperbolic equation of Kirchhoff
	type with a small parameter $\ep$ in front of the second-order
	time-derivative. 
	
	In a recent paper, under a suitable assumption on initial data, we
	proved decay-error estimates for the difference between solutions
	of the hyperbolic problem and the corresponding solutions of the
	limit parabolic problem.  These estimates show in the same time
	that the difference tends to zero both as $\ep\to 0^{+}$, and as
	$t\to +\infty$.  In particular, in that case the difference decays
	\emph{faster} than the two terms separately.
	
	In this paper we consider the complementary assumption on initial
	data, and we show that now the \emph{optimal} decay-error
	estimates involve a decay rate which is \emph{slower} than the
	decay rate of the two terms. 
	
	In both cases, the improvement or deterioration of decay rates 
	depends on the smallest frequency represented in the Fourier 
	components of initial data.
	
\vspace{1cm}

\noindent{\bf Mathematics Subject Classification 2000 (MSC2000):}
35B25, 35L70, 35L80.

\vspace{1cm} 

\end{abstract}

 
\section{Introduction}

Let $H$ be a separable real Hilbert space.  For every $x$ and $y$ in
$H$, $|x|$ denotes the norm of $x$, and $\langle x,y\rangle$ denotes
the scalar product of $x$ and $y$.  Let $A$ be a self-adjoint linear
operator on $H$ with dense domain $D(A)$.  We assume that $A$ is
nonnegative, namely $\langle Ax,x\rangle\geq 0$ for every $x\in D(A)$,
so that for every $\alpha\geq 0$ the power $A^{\alpha}x$ is defined
provided that $x$ lies in a suitable domain $D(A^{\alpha})$.

We consider the Cauchy problem
\begin{equation}
	\ep\uep''(t)+\uep'(t)+\auqg{\uep(t)}A\uep(t)=0
	\quad\quad
	\forall t\geq 0,
	\label{pbm:h-eq}
\end{equation}
\begin{equation}
	\uep(0)=u_0,\hspace{3em}\uep'(0)=u_1,
	\label{pbm:h-data}
\end{equation}
where $\ep>0$ and $\gamma\geq 1$ are real parameters, and 
$(u_{0},u_{1})\in\da\times\dau$ are initial data satisfying the 
mild nondegeneracy condition
\begin{equation}
	A^{1/2}u_{0}\neq 0.
	\label{hp:mdg}
\end{equation}

The singular perturbation problem in its generality consists in
proving the convergence of solutions of (\ref{pbm:h-eq}),
(\ref{pbm:h-data}) to solutions of the first order problem
\begin{equation}
	u'(t)+ \auqg{u(t)}Au(t)=0
	\quad\quad
	\forall t\geq 0,
	\label{pbm:p-eq}
\end{equation}
\begin{equation}
	u(0)=u_{0},
	\label{pbm:p-data}
\end{equation}
obtained setting formally $\ep=0$ in (\ref{pbm:h-eq}), and
omitting the second initial condition in~(\ref{pbm:h-data}).

This problem has generated a considerable literature in the last 30
years.  In particular, it is well known that the parabolic problem
(\ref{pbm:p-eq}), (\ref{pbm:p-data}) has a unique global solution for
every $u_{0}\in D(A)$ (and even for less regular data), and the
hyperbolic problem (\ref{pbm:h-eq}), (\ref{pbm:h-data}) has a unique
global solution provided that $(u_{0},u_{1})\in\da\times\dau$ satisfy
the nondegeneracy assumption (\ref{hp:mdg}) and $\ep$ is small enough.
The interested reader is referred to the survey~\cite{gg:survey-diss},
where more general nonlinearities and more general dissipative terms
have also been considered (see also~\cite{gg:de-dg1} for the special 
case of (\ref{pbm:h-eq}), and for more references).

Next step consists in estimating the behavior of $u(t)$, $\uep(t)$,
and of the difference $\uep(t)-u(t)$ as $t\to +\infty$ and as $\ep\to
0^{+}$. This gave rise to three types of results.
\begin{enumerate}
	\renewcommand{\labelenumi}{(\Alph{enumi})}
	\item  \emph{Decay estimates}. In this case $\ep$ is fixed, and 
	$t\to+\infty$. The prototype of decay estimates, in the case of 
	coercive operators, is that
	$$|A^{1/2}u(t)|^{2}\leq\frac{C}{(1+t)^{1/\gamma}}
	\quad\mbox{and}\quad
	|A^{1/2}\uep(t)|^{2}\leq\frac{C}{(1+t)^{1/\gamma}}$$
	for every $t\geq 0$, where the constant $C$ is independent of 
	$\ep$ and of course also of $t$. As a consequence we have also 
	that
	\begin{equation}
		|A^{1/2}(\uep(t)-u(t))|^{2}\leq\frac{C}{(1+t)^{1/\gamma}}
		\quad\quad
		\forall t\geq 0.
		\label{est:decay-proto}
	\end{equation}
	
	These estimates have been obtained for the first time
	in~\cite{ny}, and then in~\cite{gg:k-decay} in full generality.
	We point out that the decay rate of $|A^{1/2}u(t)|$ and
	$|A^{1/2}\uep(t)|$ coincides with the decay rate of solutions of
	the ordinary differential equation
	$$y'(t)+|y(t)|^{2\gamma}y(t)=0,$$
	which is actually the special case of (\ref{pbm:p-eq}) where
	$H=\re$ and $A$ is the identity (see also~\cite{haraux} for decay
	rates of second order ordinary differential equations).

	\item  \emph{Error estimates}. In this case $t$ is fixed, and 
	$\ep\to 0^{+}$. The prototype of error estimates is that for 
	initial data $(u_{0},u_{1})\in D(A^{3/2})\times D(A^{1/2})$ one 
	has that 
	\begin{equation}
		|A^{1/2}(\uep(t)-u(t))|^{2}\leq C\ep^{2}
		\quad\quad
		\forall t\geq 0,
		\label{est:error-proto}
	\end{equation}
	where the constant $C$ is once again independent of $\ep$ and $t$ 
	(global-in-time error estimates). It is well-known that $\ep^{2}$ 
	is the best possible convergence rate (even when looking for 
	local-in-time error estimates), and that $D(A^{3/2})\times 
	D(A^{1/2})$ is the minimal requirement on initial data which 
	guarantees this rate (even in the case of linear equations). We 
	refer to~\cite{gg:l-cattaneo} for these aspects.

	\item  \emph{Decay-error estimates}. They are the ultimate goal 
	of the theory, since they represent the meeting point of (A) and 
	(B). The general form of a decay-error estimate is something like
	$$|A^{1/2}(\uep(t)-u(t))|^{2}\leq\omega(\ep)\sigma(t),$$
	where of course $\omega(\ep)\to 0$ as $\ep\to 0^{+}$, and 
	$\sigma(t)\to 0$ as $t\to +\infty$.
\end{enumerate}

Trivial decay-error estimates can be obtained by exploiting the 
classical inequality $\min\{a,b\}\leq a^{1-\theta}b^{\theta}$, which 
holds true for every pair of positive real numbers $a$ and $b$, and 
every $\theta\in[0,1]$. Thus (\ref{est:decay-proto}) and 
(\ref{est:error-proto}) can be interpolated by the estimates
\begin{equation}
	|A^{1/2}(\uep(t)-u(t))|^{2}\leq 
	C\frac{\ep^{2-2\theta}}{(1+t)^{\theta/\gamma}}
	\quad\quad
	\forall t\geq 0,
	\label{de-low-cost}
\end{equation}
for a suitable constant $C$ independent of $\ep$, $t$, $\theta$. 
Nevertheless, such estimates are in general non-optimal, both with 
respect to the convergence rate $2-2\theta$, and with respect to the 
decay rate $\theta/\gamma$.

Indeed, both in the case of linear equations~\cite{ch}, and in the
case of nondegenerate Kirchhoff
equations~\cite{yamazaki,yamazaki-wd,yamazaki-cwd,gg:w-ndg}, it was
possible to prove decay-error estimates with the same decay rate of
decay estimates, and the same convergence rate of error estimates.  In
the degenerate case, this leads to guess that
\begin{equation}
	|A^{1/2}(\uep(t)-u(t))|^{2}\leq 
	C\frac{\ep^{2}}{(1+t)^{1/\gamma}}
	\quad\quad
	\forall t\geq 0.
	\label{est:de-naive}
\end{equation}

This problem was addressed for the first time in~\cite{gg:de-dg1},
where two results have been proved.  The first one is a counterexample
showing that (\ref{est:de-naive}) is false for general initial data
$(u_{0},u_{1})\in D(A^{3/2})\times D(A^{1/2})$.  The second one is
that a special assumption on initial data, which rules out the
previous counterexample, yields an estimate which is even better than
(\ref{est:de-naive}), namely that
\begin{equation}
	|A^{1/2}(\uep(t)-u(t))|^{2}\leq 
	C\frac{\ep^{2}}{(1+t)^{\delta/\gamma}}
	\quad\quad
	\forall t\geq 0
	\label{est:de-delta}
\end{equation}
for some $\delta>1$. In other words, for these data the difference 
decays faster than the two terms separately.

In order to clarify the special assumption on initial data, let us 
assume for simplicity that the spectrum of $A$ consists of a sequence 
of positive eigenvalues. Let us assume that $v_{0}$ is an eigenvector 
with respect to some eigenvalue $\mu$. Let 
$\nu>\mu$, and let us assume that
\begin{equation}
	\begin{array}{rcl}
		u_{0} & = & v_{0}+\mbox{components w.r.t.\ eigenvalues 
		$\geq\nu$,}  \\
		\noalign{\vspace{0.5ex}} u_{1} & = & \beta
		v_{0}+\mbox{components w.r.t.\ eigenvalues
		$\geq\nu$}
	\end{array}
	\label{hp:data-improve}
\end{equation}
for some $\beta\in\re$ (possibly equal to 0).  Then the second result
in~\cite{gg:de-dg1} is that (\ref{est:de-delta}) holds true with
$\delta:=\min\{2\gamma+1,\nu/\mu\}$ (with a logarithmic correction when
the two terms in the minimum are equal).  Note that $\delta>1$ in this
case.

On the contrary, the counterexample presented in~\cite{gg:de-dg1} is
that (\ref{est:de-naive}) cannot be true when $u_{0}$ and $u_{1}$ are
eigenvectors of $A$ with respect to eigenvalues $\mu$ and $\nu$,
respectively, with $\nu<\mu$.  In other words, in that case a
deterioration of decay rates is expected.  Quantifying this
deterioration was stated as Open Problem~2 in section~4
of~\cite{gg:de-dg1}.

In this paper we investigate this open problem, and we find the best
decay-error estimates which are true in this case.  So we assume that
$v_{0}$ is an eigenvector of $A$ with respect to some eigenvalue
$\mu$, and that $v_{1}$ is an eigenvector of $A$ with respect to some
eigenvalue $\nu$, with $\nu\leq\mu$.  Then we take initial data of the
form
\begin{equation}
	\begin{array}{rcl}
		u_{0} & = & v_{0}+\mbox{components w.r.t.\ eigenvalues 
		$>\mu$,}  \\
		\noalign{\vspace{0.5ex}}
		u_{1} & = & v_{1}+\mbox{components w.r.t.\ eigenvalues
		$>\nu$.}
	\end{array}
	\label{hp:data-here}
\end{equation}

In this case we prove that (\ref{est:de-delta}) holds true with
$\delta:=\nu/\mu$, which is now less than or equal to 1.  We also show
that this exponent is optimal, in the sense that we cannot obtain a
better decay rate if we want to keep the optimal convergence rate
$\ep^{2}$.

When the spectrum of $A$ is an increasing sequence of eigenvalues
(with no assumption on the dimension of eigenspaces),
assumption~(\ref{hp:data-here}) is complementary to
assumption~(\ref{hp:data-improve}).  This happens, for example, when
$A$ is the Dirichlet-Laplacian in a bounded open set with reasonable
boundary, namely the operator involved in Kirchhoff equations in
concrete form.  Therefore, for these operators the problem of decay-error
estimates is now fully closed by Theorem~2.5 in~\cite{gg:de-dg1}, and
Theorem~\ref{thm:main} of the present paper.  We refer to
Remark~\ref{rmk:data} for the details.

This paper is organized as follows.  In section~\ref{sec:main} we fix
the notation and state our main results.  In
section~\ref{sec:heuristic} we present a toy model which roughly
explains the deterioration of decay rates.  In
section~\ref{sec:proofs} we prove our results.

\setcounter{equation}{0}

\section{Statements}\label{sec:statements}

\subsection{Notation and statements}\label{sec:main}

In the following we always assume that $H$ is a Hilbert space, and $A$
is a nonnegative self-adjoint (unbounded) operator on $H$ with dense
domain.  We always assume that $\gamma\geq 1$ is a real number, and
that $(u_{0},u_{1})\in\da\times\dau$ is a pair of initial conditions
satisfying the nondegeneracy condition (\ref{hp:mdg}).

Let $u(t)$ be the unique global solution of the first order problem
(\ref{pbm:p-eq}), (\ref{pbm:p-data}).  Let $\uep(t)$ be the unique
global solution of the second order problem (\ref{pbm:h-eq}),
(\ref{pbm:h-data}), which exists at least for every
$\ep\in(0,\ep_{0})$ for some positive $\ep_{0}$.  Following the
approach introduced in~\cite{lions} in the linear case, we define the
corrector $\tetep(t)$ as the solution of the second order
\emph{linear} ordinary differential equation
\begin{equation}
	\ep\tetep''(t)+\tetep'(t)=0 \hspace{2em}
	\forall t\geq 0,
	\label{pbm:tetep-eq}
\end{equation}
with initial data
\begin{equation}
	\tetep(0)=0,\hspace{2em}\tetep'(0)=u_1+
	\auqg{u_{0}}Au_{0}=:\theta_{0}.
	\label{pbm:tetep-data}
\end{equation}

Since $\theta_{0}=\uep'(0)-u'(0)$, this corrector keeps into account
the boundary layer due to the loss of one initial condition.  

Finally, we define $\rep(t)$ and $\roep(t)$ in such a way that
$$\uep(t)=u(t)+\tetep(t)+\rep(t)=u(t)+\roep(t)\quad\quad\forall t\geq
0.$$

With these notations, the singular perturbation problem consists in
proving that $\rep(t)\to 0$ or $\roep(t)\to 0$ in some sense as
$\ep\to 0^{+}$.  We recall that the two remainders play different
roles.  In particular, $\rep(t)$ is well suited for estimating
derivatives, while $\roep(t)$ is used in estimates without
derivatives.  This distinction is essential.  Indeed it is not
possible to prove decay-error estimates on $A^{\alpha}\rep(t)$ because
it does not decay to 0 as $t\to +\infty$ (indeed $\uep(t)$ and $u(t)$
tend to 0, while the corrector $\tetep(t)$ does not), and it is not
possible to prove decay-error estimates on $A^{\alpha}\roep'(t)$
because in general for $t=0$ it does not tend to 0 as $\ep\to 0^{+}$
(due to the loss of one initial condition).

Now let us introduce our assumptions on initial data.  To this end, we
need some basic facts from the spectral theory of operators, which we
recall following~\cite{rudin}.

Let $E$ be the resolution of the identity associated with the operator
$A$.  For every measurable subset $J\subseteq[0,+\infty)$ we consider
the space $H_{J}:=\mathcal{R}(E(J))$, namely the range of the
projection operator $E(J)$, which is a closed subspace of $H$.  In the
case where $H$ admits a (finite or countable) orthonormal system
$\{e_{k}\}$ made by eigenvalues of $A$, and $\{\lambda_{k}^{2}\}$ is
the sequence of corresponding eigenvalues, then $H_{J}$ is just the
set of all $v\in H$ such that $\langle v,e_{k}\rangle=0$ for every
$k\in\n$ such that $\lambda_{k}^{2}\not\in J$.

We are now ready to introduce the class of initial data we consider in
this paper.

\begin{defn}[Assumption on initial data]\label{hp:data}
	\begin{em}
		Let $\mu$ and $\nu$ be two positive real numbers.  We say that
		a pair of initial conditions $(u_{0},u_{1})\in\da\times\dau$
		satisfies the $(\mu,\nu)$-assumption if we can write 
		$u_{0}=v_{0}+w_{0}$ and $u_{1}=v_{1}+w_{1}$, where
		\begin{itemize}
			\item $v_{0}$ and $v_{1}$ are eigenvectors of $A$ with 
			respect to the eigenvalues $\mu$ and $\nu$, respectively,
			
			\item $w_{0}\in H_{(\mu,+\infty)}$ and $w_{1}\in
			H_{(\nu,+\infty)}$.
		\end{itemize}
	\end{em}
\end{defn}

In other words, $\mu$ is the smallest frequency with respect to which
$u_{0}$ has a nonzero component, and $v_{0}$ is such a component.
Analogously, $\nu$ is the smallest frequency with respect to which
$u_{1}$ has a nonzero component $v_{1}$.

The main result of this paper is the following.

\begin{thm}[Decay-error estimates]
	\label{thm:main}
	Let $H$, $A$, $\gamma$, $(u_{0},u_{1})$, $\ep_{0}$, $u(t)$,
	$\uep(t)$, $\tetep(t)$, $\roep(t)$, $\rep(t)$ be as usual. 
	
	Let us assume that the pair $(u_{0},u_{1})$ satisfies the
	$(\mu,\nu)$-assumption with
	$$0<\nu\leq\mu.$$
	
	Then the following conclusions hold true with $\delta:=\nu/\mu$.
	\begin{enumerate}
		\renewcommand{\labelenumi}{(\arabic{enumi})}
		\item If in addition we assume that
		$(u_{0},u_{1})\in\dat\times\dau$, then there exist
		$\ep_{1}\in(0,\ep_{0})$ and a constant $C$ such that for
		every $\ep\in(0,\ep_{1})$ we have that
		$$|\roep(t)|^{2}+|A^{1/2}\roep(t)|^{2}+ \ep(1+t)|\rep'(t)|^{2}
		\leq C\frac{\ep^{2}}{(1+t)^{\delta/\gamma}}
		\quad\quad
		\forall t\geq 0,$$
		$$\int_{0}^{t}(1+s)^{2\delta/\gamma}\left((1+s)|\rep'(s)|^{2}+
		\frac{|A^{1/2}\roep(s)|^{2}}{1+s}\right)\,ds\leq 
		C\ep^{2}(1+t)^{\delta/\gamma}
		\quad\quad
		\forall t\geq 0.$$
		
		\item If in addition we assume that $(u_{0},u_{1})\in
		D(A^{2})\times\da$, then there exist $\ep_{1}\in(0,\ep_{0})$
		and a constant $C$ such that for every $\ep\in(0,\ep_{1})$ we
		have that
		$$|A\roep(t)|^{2}+(1+t)^{2}|\rep'(t)|^{2}
		\leq C\frac{\ep^{2}}{(1+t)^{\delta/\gamma}}
		\quad\quad
		\forall t\geq 0,$$
		$$\int_{0}^{t}(1+s)^{2\delta/\gamma}\left((1+s)|A^{1/2}\rep'(s)|^{2}+
		\frac{|A\roep(s)|^{2}}{1+s}\right)\,ds\leq 
		C\ep^{2}(1+t)^{\delta/\gamma}
		\quad\quad
		\forall t\geq 0.$$
	\end{enumerate}
\end{thm}

As we are going to see in the proofs, the main point is obtaining the
estimate on the lowest order term $\roep(t)$.  Quite surprisingly, all
remaining estimates follow from this one through a linear argument
presented in~\cite{gg:de-dg1}.

The following result shows that the estimates provided by 
Theorem~\ref{thm:main} (at least the one on $\roep(t)$, which implies 
all the rest) are optimal.

\begin{thm}[Optimality of decay-error estimates]
	\label{thm:optimality}
	Let $H$, $A$, $\gamma$, $(u_{0},u_{1})$, $\ep_{0}$, $u(t)$,
	$\uep(t)$, $\tetep(t)$, $\roep(t)$, $\rep(t)$ be as usual. 
	
	Let us assume that $u_{0}$ and $u_{1}$ are eigenvectors of $A$
	with respect to the eigenvalues $\mu$ and $\nu$, respectively, and
	let $\delta:=\nu/\mu$.  Let us assume also that
	\begin{itemize}
		\item  either $\nu<\mu$,
	
		\item  or $\nu=\mu$, and $u_{0}$ and $u_{1}$ are orthogonal.
	\end{itemize}
	
	Then there exist $\ep_{1}\in(0,\ep_{0})$ and $C>0$ such that for
	every $\ep\in(0,\ep_{1})$ we have that
	\begin{equation}
		\sup_{t\geq 0}\left\{(1+t)^{\delta/\gamma}|\roep(t)|^{2}\right\}\geq 
		C\ep^{2}.
		\label{th:optimality}
	\end{equation}
\end{thm}

We conclude with some comments about possible extensions and 
applications.

\begin{rmk}
	\begin{em}
		For the sake of simplicity, we limit ourselves to state and
		prove (\ref{th:optimality}) when initial data are eigenvectors.
		The result is actually true for all initial data of the
		form~(\ref{hp:data-here}).  The reason is that components
		corresponding to higher frequencies decay faster
		(see~\cite{ghisi:decay}), and therefore the decay rate is
		always dictated by the smallest frequencies represented in
		$u_{0}$ and $u_{1}$.
		
		Similarly, in the case where $\nu=\mu$, we could weaken the 
		assumption that $u_{0}$ and $u_{1}$ are orthogonal by just 
		asking that they are linearly independent. In this case the 
		proof should be modified in an obvious way by introducing 
		the component of $u_{1}$ orthogonal to $u_{0}$.
	\end{em}
\end{rmk}

\begin{rmk}\label{rmk:data}
	\begin{em}
		Let us assume that the spectrum of $A$ consists of an
		increasing sequence of positive eigenvalues.  Let
		$(u_{0},u_{1})\in D(A)\times D(A^{1/2})$ be any pair of
		initial conditions satisfying the nondegeneracy condition
		(\ref{hp:mdg}).  Let $\mu$ and $\nu$ be the smallest
		eigenvalues with respect to which $u_{0}$ and $u_{1}$ have
		nonzero components $v_{0}$ and $v_{1}$, respectively.
		\begin{itemize}
			\item If $\nu>\mu$ (or if $u_{1}=0$, in which case $\nu$ 
			is not well defined) we are in the situation
			of~\cite{gg:de-dg1}, which yields decay-error estimates
			with improved decay rates.
		
			\item If $\nu<\mu$ we are in the assumptions of
			Theorem~\ref{thm:main} above, which yields decay error
			estimates with deteriorated decay rates, and in general
			nothing more because of Theorem~\ref{thm:optimality}
			above.  
		
			\item If $\nu=\mu$ we are for sure in the situation of
			Theorem~\ref{thm:main}, which in any case guarantees
			decay-error estimates without improvement or deterioration
			of decay rates ($\delta=1$).  If $v_{0}$ and $v_{1}$ are
			linearly independent, then $\delta=1$ is optimal because
			of Theorem~\ref{thm:optimality}.  If $v_{0}$ and $v_{1}$
			are linearly dependent, then we are also in the
			assumptions of~\cite{gg:de-dg1}, and once again we get an
			improvement of decay rates.
		\end{itemize}
		
	\end{em}
\end{rmk}

\begin{rmk}
	\begin{em}
		The previous remark extends easily to operators whose 
		spectrum is an increasing sequence of nonnegative eigenvalues 
		\emph{including} 0 (this is the case, for example, of the 
		Neumann-Laplacian on a bounded interval).
		
		Indeed in this case it is enough to separate components with 
		respect to the kernel, where $u(t)$ is constant and 
		$\uep(t)$ coincides with the corrector $\tetep(t)$, and apply 
		the theory of this paper in the subspace orthogonal to the 
		kernel. We point out that, since the nonlinear term depends 
		on $A^{1/2}u$ or $A^{1/2}\uep$, what happens in the kernel 
		has no influence in the orthogonal subspace.
	\end{em}
\end{rmk}

\subsection{Heuristics}\label{sec:heuristic}

In this section we present a simple calculation on ordinary
differential equations that leads to guess that
$|\roep(t)|\sim\ep(1+t)^{-\delta/(2\gamma)}$.  Let us start with two
simplifications.
\begin{itemize}
	\item  When $\ep$ is small enough, the parabolic equation is a 
	good approximation of the hyperbolic one. As a consequence, after 
	the initial layer, $\uep(t)$ and $u(t)$ can \emph{both} be 
	considered as solutions of the \emph{parabolic} equation.

	\item After the initial layer, for example at time $t=1$, the
	difference $\uep(t)-u(t)$ is of order $\ep$.  Actually this
	follows from the well-established local-in-time error estimates.
\end{itemize}

Under these simplifying assumptions, the singular perturbation 
problem is reduced to estimating the difference between two solutions 
of the \emph{parabolic} equation whose initial data differ of order 
$\ep$.

So let $v_{0}$ and $v_{1}$ be orthonormal eigenvectors of $A$ with
respect to eigenvalues $\mu$ and $\nu$, respectively, with
$0<\nu<\mu$.  Let us consider the solutions $u_{0}(t)$ and $u_{1}(t)$
of the parabolic problem with initial data $u_{0}(0)=v_{0}$, and
$u_{1}(0)=v_{0}+\ep v_{1}$, respectively (so that at time $t=0$ the
difference is of order $\ep$, and lies on a component corresponding to
the smallest frequency).

Now it is easy to see that $u_{0}(t)$ is always a multiple of 
$v_{0}$, while $u_{1}(t)$ can be written in the form 
$u_{1}(t):=w(t)v_{0}+v(t)v_{1}$, where $v(t)$ and $w(t)$ are 
solutions of the system of ordinary differential equations
$$\left\{
\begin{array}{l}
	w'(t)+\mu\left(\strut \nu v^{2}(t)+\mu
	w^{2}(t)\right)^{\gamma}w(t)=0, \\
	\noalign{\vspace{1ex}}
	v'(t)+\nu\left(\strut \nu v^{2}(t)+\mu
	w^{2}(t)\right)^{\gamma}v(t)=0,
\end{array}
\right.$$
with initial data $w(0)=1$, and $v(0)=\ep$.

Now we make a further simplifying assumption, namely that 
$u_{1}(t)-u_{0}(t)\sim v(t)v_{1}$, which is reasonable if we accept 
that components corresponding to lower frequencies decay more slowly. 
Thus we have reduced ourselves to estimating $v(t)$.

According to the main trick introduced in~\cite{k-par}, $v(t)$ and
$w(t)$ can be written in the form 
$$v(t):=\ep\psiep(t), \hspace{3em}
w(t):=\left[\psiep(t)\right]^{1/\delta},$$
where as usual $\delta:=\nu/\mu$, and $\psiep(t)$ solves
\begin{equation}
	\psiep'(t)+\nu\left(\nu\ep^{2}\psiep^{2}(t)+
	\mu\psiep^{2/\delta}(t)\right)^{\gamma}\psiep(t)=0,
	\quad\quad
	\psiep(0)=1.
	\label{eqn:psi}
\end{equation}

Now we make the final simplifying assumption by setting $\ep=0$ in
(\ref{eqn:psi}), so that from now on $\psiep(t)$ does not depend on
$\ep$ and solves
$$\psi'(t)+k\psi^{1+2\gamma/\delta}(t)=0,
\quad\quad
\psi(0)=1,$$
for a suitable positive constant $k$.  This differential equation can
be easily integrated, giving that
$\psi(t)\sim(1+t)^{-\delta/(2\gamma)}$, hence
$$|\roep(t)|\sim|u_{1}(t)-u_{0}(t)|\sim
|v(t)|=\ep\psi(t)\sim \ep(1+t)^{-\delta/(2\gamma)},$$
which is consistent both with Theorem~\ref{thm:main}, and with 
Theorem~\ref{thm:optimality}.

We conclude with some remarks. First of all, (\ref{eqn:psi}) is more 
or less the same equation which appears in the proof of 
Theorem~\ref{thm:optimality} (see also Lemma~\ref{lemma:ODE}), which 
means that the rough calculation we did is actually close to reality.

Secondly, neglecting the term with $\ep$ in (\ref{eqn:psi}) is not 
reasonable for all times. Indeed $\psiep(t)\to 0$ as $t\to 
+\infty$, and in this regime the term with $\psiep^{2}$ becomes 
dominant over $\psiep^{2/\delta}$ (because $\delta<1$). This means that 
as $t\to+\infty$ the true approximation of (\ref{eqn:psi}) is 
$$\psiep'(t)+k\ep^{2\gamma}\psiep^{1+2\gamma}(t)=0
\quad\quad
\psiep(0)=1,$$
which gives $\psi(t)\sim(1+\ep^{2\gamma}t)^{-1/(2\gamma)}$, namely a
better decay rate.

Which is the correct approximation of (\ref{eqn:psi})?  In a certain
sense both!  The fact that
$|\roep(t)|\sim\ep(1+\ep^{2\gamma}t)^{-1/(2\gamma)}$, suggested by the
second approximation, is consistent with previous works where it was
proved that both $u(t)$ and $\uep(t)$ decay as $(1+t)^{-1/(2\gamma)}$.
The point is that the final coefficient in this case turns out to be
of order 1, not of order $\ep$, and therefore we get a better decay
rate which we pay by losing the convergence rate.  The first
approximation, on the contrary, preserves the optimal convergence rate
of order $\ep$. As in (\ref{de-low-cost}), one could interpolate the 
two extremes by finding a family of estimates with intermediate decay 
rates and intermediate convergence rates.

The interplay between the two different regimes of (\ref{eqn:psi}) is
what makes highly nontrivial the rigorous analysis of decay-error 
estimates for nonlinear degenerate equations.

\setcounter{equation}{0}
\section{Proofs}\label{sec:proofs}

In the sequel we set
\begin{equation}
	c(t):=|A^{1/2}u(t)|^{2\gamma},
	\hspace{3em}
	\cep(t):=|A^{1/2}\uep(t)|^{2\gamma},
	\label{defn:c-cep}
\end{equation}
and we consider their anti-derivatives
\begin{equation}
	C(t):=\int_{0}^{t}c(s)\,ds,
	\hspace{3em}
	\Cep(t):=\int_{0}^{t}\cep(s)\,ds.
	\label{defn:C-Cep}
\end{equation}

We exploit that the corrector $\tetep(t)$, which is the solution of 
(\ref{pbm:tetep-eq}), (\ref{pbm:tetep-data}), is given by the explicit 
formula
\begin{equation}
	\tetep(t)=\ep \theta_{0}
	\left(1-e^{-t/\ep}\strut\right)
	\quad\quad
	\forall t\geq 0.
	\label{eqn:tetep}
\end{equation}

We also set
\begin{equation}
	\gep(t):=-(\cep(t)-c(t))Au(t)-\ep u''(t),
	\label{defn:gep}
\end{equation}
so that we can regard $\roep(t)$ as the solution of the linear equation
\begin{equation}
	\ep\roep''(t)+\roep'(t)+\cep(t)A\roep(t)=\gep(t),
	\label{eqn:roep-lin}
\end{equation}
with initial data
\begin{equation}
	\roep(0)=0,
	\quad\quad
	\roep'(0)=\theta_{0},
	\label{eqn:roep-data}
\end{equation}
and $\rep(t)$ as the solution of the linear equation
\begin{equation}
	\ep\rep''(t)+\rep'(t)+\cep(t)A\roep(t)=\gep(t),
	\label{eqn:rep-lin}
\end{equation}
with initial data
\begin{equation}
	\rep(0)=0,
	\quad\quad
	\rep'(0)=0.
	\label{eqn:rep-data}
\end{equation}

In many points we need to consider components of $\uep(t)$ or $u(t)$ 
with respect to $v_{0}$ and $v_{1}$, where $v_{0}$ and $v_{1}$ are 
the eigenvectors of $A$ which appear in the $(\mu,\nu)$-assumption on 
initial data. To this end we set
\begin{equation}
	\uepn(t):=\frac{1}{|v_{1}|}\langle\uep(t),v_{1}\rangle,
	\hspace{3em}
	\uepm(t):=\frac{1}{|v_{0}|}\langle\uep(t),v_{0}\rangle.
	\label{defn:uepn-uepm}
\end{equation}

It is easy to see that these real functions satisfy, respectively, the
following second order equations
\begin{equation}
	\ep\uepn''(t)+\uepn'(t)+\nu\cep(t)\uepn(t)=0,
	\quad\quad
	\ep\uepm''(t)+\uepm'(t)+\mu\cep(t)\uepm(t)=0.
	\label{eqn:uepn-uepm}
\end{equation}

\subsection{Previous results}

Many estimates on $u(t)$ and $\uep(t)$ have been proved in literature.
We refer to~\cite{gg:k-decay} (see also Theorems~A and~B
in~\cite{gg:de-dg1}) for a complete list.  For the convenience of the
reader, in next result we limit ourselves to recall just those
estimates needed in the sequel.

\begin{thmbibl}\label{bibl:main}
	Let $H$, $A$, $\gamma$, $(u_{0},u_{1})$, $\ep_{0}$, $u(t)$,
	$\uep(t)$, $\tetep(t)$, $\roep(t)$, $\rep(t)$ be as usual.  Let us
	assume that $u_{0}$ and $u_{1}$ are in $H_{[\nu,+\infty)}$ for
	some $\nu>0$ (which is equivalent to asking the coercivity of $A$).
	
	Then there exists $\ep_{1}>0$, and positive constants 
	$M_{1},\ldots,M_{16}$, such that the following estimates 
	hold true for every $\ep\in(0,\ep_{1})$.
	\begin{enumerate}
		\renewcommand{\labelenumi}{(\arabic{enumi})}
		\item \emph{(Decay estimates for the parabolic problem)} We
		have that
		\begin{equation}
			|A^{1/2}u(t)|^{2}\leq\frac{M_{1}}{(1+t)^{1/\gamma}}
			\quad\quad
			\forall t\geq 0,
			\label{th:Aqu}
		\end{equation}
		\begin{equation}
			|Au(t)|^{2}\leq\frac{M_{2}}{(1+t)^{1/\gamma}}
			\quad\quad
			\forall t\geq 0,
			\label{th:Au}
		\end{equation}
		and as a consequence
		\begin{equation}
			c(t)\leq\frac{M_{3}}{1+t}
			\quad\quad
			\forall t\geq 0.
			\label{th:c}
		\end{equation}
		
		If in addition $u_{0}\in D(A^{3/2})$, then for every 
		$\delta\in(0,2\gamma+1]$ we have that
		\begin{equation}
			|A^{3/2}u(t)|^{2}\leq\frac{M_{4}}{(1+t)^{1/\gamma}}
			\quad\quad
			\forall t\geq 0,
			\label{th:Atu}
		\end{equation}
		\begin{equation}
			\int_{0}^{t}(1+s)^{1+2\delta/\gamma}|u''(s)|^{2}\,ds\leq
			M_{5}(1+t)^{\delta/\gamma}
			\quad\quad
			\forall t\geq 0.
			\label{th:u''-int}
		\end{equation}
	
		If in addition $u_{0}\in D(A^{2})$, then for every 
		$\delta\in(0,2\gamma+1]$ we have that 
		\begin{equation}
			\int_{0}^{t}(1+s)^{1+2\delta/\gamma}|A^{1/2}u''(s)|^{2}\,ds\leq
			M_{6}(1+t)^{\delta/\gamma}
			\quad\quad
			\forall t\geq 0,
			\label{th:auq-u''-int}
		\end{equation}
		\begin{equation}
			|u''(t)|^{2}\leq\frac{M_{7}}{(1+t)^{4+1/\gamma}}
			\quad\quad
			\forall t\geq 0.
			\label{th:u''}
		\end{equation}
	
		\item \emph{(Decay estimates for the hyperbolic problem)} We
		have that
		\begin{equation}
			\frac{M_{8}}{(1+t)^{1/\gamma}}\leq
			|A^{1/2}\uep(t)|^{2}\leq
			\frac{M_{9}}{(1+t)^{1/\gamma}}
			\quad\quad
			\forall t\geq 0,
			\label{th:Aquep}
		\end{equation}
		\begin{equation}
			|A\uep(t)|^{2}\leq
			\frac{M_{10}}{(1+t)^{1/\gamma}}
			\quad\quad
			\forall t\geq 0,
			\label{th:Auep}
		\end{equation}
		\begin{equation}
			|\uep'(t)|^{2}\leq\frac{M_{11}}{(1+t)^{2+1/\gamma}}
			\quad\quad
			\forall t\geq 0,
			\label{th:uep'}
		\end{equation}
		and as a consequence
		\begin{equation}
			\frac{M_{12}}{1+t}\leq\cep(t)\leq\frac{M_{13}}{1+t}
			\quad\quad
			\forall t\geq 0,
			\label{th:cep}
		\end{equation}
		\begin{equation}
			\frac{|\cep'(t)|}{\cep(t)}\leq\frac{M_{14}}{1+t}
			\quad\quad
			\forall t\geq 0.
			\label{th:cep'}
		\end{equation}
	
		\item  \emph{(Basic error estimates)} If in addition 
		$(u_{0},u_{1})\in D(A^{3/2})\times D(A^{1/2})$, we have that
		\begin{equation}
			|\roep(t)|^{2}\leq M_{15}\ep^{2}
			\quad\quad
			\forall t\geq 0,
			\label{th:ee-roep}
		\end{equation}
		\begin{equation}
			\int_{0}^{+\infty}(1+t)|\rep'(t)|^{2}\,dt\leq M_{16}\ep^{2}.
			\label{th:ee-rep'}
		\end{equation}
	
		\item  \emph{(``Monotonicity'' estimates)} We have that
		\begin{equation}
			\langle\cep(t)A\uep(t)-c(t)Au(t),\roep(t)\rangle\geq
			\frac{1}{2}\left(\cep(t)+c(t)\right)|A^{1/2}\roep(t)|^{2}
			\quad\quad
			\forall t\geq 0.
			\label{th:monoton}
		\end{equation}
	\end{enumerate}
	
\end{thmbibl}

We refer to \cite{gg:k-decay} and to Lemma~3.3 in~\cite{gg:de-dg1} for
the estimates on solutions of the parabolic and hyperbolic problem, to
Theorem~2.4 in~\cite{ghisi:error} (see also~\cite{gg:k-PS}) for the
basic error estimates, and to Lemma~3.4 in~\cite{ghisi:error} for the
monotonicity estimate (which actually is a general property of
vectors).  We point out that eigenvalues and components play no role
in Theorem~\ref{bibl:main} above.

The situation is different in the
following result, where $\uep(t)$ and its components $\uepn(t)$ and
$\uepm(t)$ are estimated in terms of $\Cep(t)$.  For a proof, we refer
to Theorem~3.1 in~\cite{ghisi:decay}, where similar estimates have
been obtained more generally for the components of $\uep(t)$ in the
subspaces $H_{[\nu,+\infty)}$ and $H_{[\mu,+\infty)}$.

\begin{thmbibl}
	Let $H$, $A$, $\gamma$, $(u_{0},u_{1})$, $\ep_{0}$, $\uep(t)$ be
	as usual.  Let us assume that the pair $(u_{0},u_{1})$ satisfies
	the $(\mu,\nu)$-assumption with respect to some $0<\nu\leq\mu$,
	let $\uepn(t)$ and $\uepm(t)$ be the components of $\uep(t)$,
	defined according to (\ref{defn:uepn-uepm}), and let $\Cep(t)$ be 
	defined according to (\ref{defn:C-Cep}).
	
	Then there exists $\ep_{1}>0$, and constants
	$M_{17}$, $M_{18}$, $M_{19}$ such that for every
	$\ep\in(0,\ep_{1})$ we have that
	\begin{equation}
		|A^{1/2}\uep(t)|^{2}\leq
		M_{17}e^{-2\nu\Cep(t)}
		\quad\quad
		\forall t\geq 0,
		\label{th:uep-Cep}
	\end{equation}
	\begin{equation}
		|\uepn(t)|^{2}+\frac{|\uepn'(t)|^{2}}{[\cep(t)]^{2}}\leq
		M_{18}e^{-2\nu\Cep(t)}
		\quad\quad
		\forall t\geq 0,
		\label{th:uepn}
	\end{equation}
	\begin{equation}
		|\uepm(t)|^{2}+\frac{|\uepm'(t)|^{2}}{[\cep(t)]^{2}}\leq
		M_{19}e^{-2\mu\Cep(t)}
		\quad\quad
		\forall t\geq 0.
		\label{th:uepm}
	\end{equation}
	
\end{thmbibl}

The last result we need allows to deduce all the conclusions of
Theorem~\ref{thm:main} from the only estimate on $\roep(t)$.  It is
actually a result on linear equations, in the sense that now we regard
$\roep(t)$ and $\rep(t)$ as solutions of (\ref{eqn:roep-lin}) and
(\ref{eqn:rep-lin}), forgetting that $c(t)$, $\cep(t)$, and $\gep(t)$
are given by (\ref{defn:c-cep}) and (\ref{defn:gep}), respectively.
For a proof, we refer to Proposition~3.5 in~\cite{gg:de-dg1}.

\begin{propbibl}\label{prop:linear}
	Let $H$ and $A$ be as usual, and let 
	$\ep_{0}$, $\gamma$, $\delta$ be positive real numbers.  
	
	For every $\ep\in (0,\ep_{0})$, let us assume that $\roep(t)$,
	$\rep(t)$, $\cep(t)$, and $\gep(t)$ satisfy (\ref{eqn:roep-lin})
	through (\ref{eqn:rep-data}).  Moreover, let us assume that
	\begin{enumerate}
		\renewcommand{\labelenumi}{(\roman{enumi})} 
		\item  the solution $\roep(t)$ satisfies the a priori estimate
		\begin{equation}
			|\roep(t)|^{2}\leq M_{20}\frac{\ep^{2}}{(1+t)^{\delta/\gamma}}
			\quad\quad
			\forall t\geq 0,
			\label{th:roep}
		\end{equation}
	
		\item the coefficient $\cep:[0,+\infty)\to(0,+\infty)$ is of
		class $C^{1}$ and satisfies (\ref{th:cep}) and 
		(\ref{th:cep'}),
	
		\item the forcing term $\gep:[0,+\infty)\to H$ is continuous
		and such that
		\begin{equation}
			\int_{0}^{t}(1+s)^{1+2\delta/\gamma}|\gep(s)|^{2}\,ds\leq
			M_{21}\ep^{2}(1+t)^{\delta/\gamma} 
			\quad\quad
			\forall t\geq 0.
			\label{hp:gep-1}
		\end{equation}
	\end{enumerate}
	
	Then the following conclusions hold true.
	\begin{enumerate}
		\renewcommand{\labelenumi}{(\arabic{enumi})}
		\item  If $\theta_{0}\in\dau$, then all the estimates in 
		statement~(1) of Theorem~\ref{thm:main} hold true.
	
		\item If in addition we have that $\theta_{0}\in\da$, and
		$\gep(t)$ satisfies also
		\begin{equation}
			\int_{0}^{t}(1+s)^{1+2\delta/\gamma}|A^{1/2}\gep(s)|^{2}\,ds\leq
			M_{22}\ep^{2}(1+t)^{\delta/\gamma}
			\quad\quad
			\forall t\geq 0,
			\label{hp:gep-2}
		\end{equation}
		\begin{equation}
			|\gep(t)|^{2}\leq
			M_{23}\frac{\ep^{2}}{(1+t)^{2+\delta/\gamma}}
			\quad\quad
			\forall t\geq 0,
			\label{hp:gep-3}
		\end{equation}
		then all the estimates in statement~(2) of
		Theorem~\ref{thm:main} hold true.
	\end{enumerate}
\end{propbibl}

\subsection{Preliminary estimates}

In the following result, we collect some further estimates on
solutions of the hyperbolic problem.  Most of them have been already
used somewhere in previous papers, but for the convenience of the
reader we give here a self contained proof (of course based on the
estimates of the previous section).  

\begin{prop}
	Let us consider the same assumptions of Theorem~\ref{thm:main},
	and let $C(t)$ and $\Cep(t)$ be defined by (\ref{defn:C-Cep}).
	
	Then there exists $\ep_{1}\in(0,\ep_{0})$, and positive constants
	$M_{24}$, \ldots, $M_{27}$, such that for every
	$\ep\in(0,\ep_{1})$ we have that
	\begin{eqnarray}
		 & e^{2\mu\gamma C(t)}\leq M_{24}(1+t)
		\quad\quad
		\forall t\geq 0, & 
		\label{th:C-above}  \\
		\noalign{\vspace{1ex}}
		 & e^{2\nu\gamma\Cep(t)}\leq M_{25}(1+t)
		\quad\quad
		\forall t\geq 0, & 
		\label{th:Cep-above}  \\
		\noalign{\vspace{1ex}}
		 & e^{2\mu\gamma C(t)}\geq M_{26}(1+t)
		\quad\quad
		\forall t\geq 0, & 
		\label{th:C-below}  \\
		\noalign{\vspace{1ex}}
		 & e^{2\mu\gamma\Cep(t)}\geq M_{27}(1+t)
		\quad\quad
		\forall t\geq 0. & 
		\label{th:Cep-below}
	\end{eqnarray}
\end{prop}

\paragraph{\emph{\textmd{Proof}}}

Let us prove the four estimates separately.

\subparagraph{\emph{\textmd{Proof of estimate (\ref{th:C-above})}}}

From (\ref{pbm:p-eq}) we have that
\begin{equation}
	\left[|A^{1/2}u(t)|^{2}e^{2\mu C(t)}\right]'=
	2\mu c(t)e^{2\mu C(t)}|A^{1/2}u(t)|^{2}-
	2e^{2\mu C(t)}c(t)|Au(t)|^{2}.
	\label{est:deriv}
\end{equation}

On the other hand, $u(t)\in H_{[\mu,+\infty)}$ for every $t\geq 0$, 
hence $|Au(t)|^{2}\geq\mu|A^{1/2}u(t)|^{2}$. It follows that the 
right-hand side of (\ref{est:deriv}) is less than or equal to 0, hence
$$|A^{1/2}u(t)|^{2}e^{2\mu C(t)}\leq|A^{1/2}u_{0}|^{2}.$$

Therefore we have that
$$\left[e^{2\mu\gamma C(t)}\right]'=
2\mu\gamma c(t)e^{2\mu\gamma C(t)}=2\mu\gamma
\left(|A^{1/2}u(t)|^{2}e^{2\mu C(t)}\right)^{\gamma}\leq
k_{1},$$
so that (\ref{th:C-above}) follows by integration.

\subparagraph{\emph{\textmd{Proof of estimate (\ref{th:Cep-above})}}}

Thanks to (\ref{th:uep-Cep}) we have that 
$$\left[e^{2\nu\gamma\Cep(t)}\right]'=
2\nu\gamma\cep(t)e^{2\nu\gamma\Cep(t)}=
2\nu\gamma\left(|A^{1/2}\uep(t)|^{2}e^{2\nu\Cep(t)}\right)^{\gamma}\leq 
k_{2},$$
so that (\ref{th:Cep-above}) follows by integration.

\subparagraph{\emph{\textmd{Proof of estimate (\ref{th:C-below})}}}

Let us begin by showing that
\begin{equation}
	|A^{1/2}u(t)|^{2}e^{2\mu C(t)}\geq k_{3}>0
	\quad\quad
	\forall t\geq 0.
	\label{est:Au-C}
\end{equation}

To this end, let $v_{0}$ denote the component of $u_{0}$ with respect
to the eigenspace of $\mu$ (as in Definition~\ref{hp:data}), and let
$u_{\mu}(t):=|v_{0}|^{-1}\langle u(t),v_{0}\rangle$ be the component
of $u(t)$ with respect to the same eigenspace.  Then of course
$u_{\mu}(t)$ satisfies 
$$u_{\mu}'(t)+\mu c(t)u_{\mu}(t)=0, \quad\quad
u_{\mu}(0)=|v_{0}|,$$
hence $u_{\mu}(t)=|v_{0}|e^{-\mu C(t)}$ for every $t\geq 0$. 
Therefore we have that
$$|A^{1/2}u(t)|^{2}e^{2\mu C(t)}\geq\mu
|u_{\mu}(t)|^{2}e^{2\mu C(t)}\geq\mu|v_{0}|^{2},$$
which implies (\ref{est:Au-C}). It follows that
$$\left[e^{2\mu\gamma C(t)}\right]'=2\mu\gamma
c(t)e^{2\mu\gamma C(t)}=2\mu\gamma
\left(|A^{1/2}u(t)|^{2}e^{2\mu C(t)}\right)^{\gamma}\geq
k_{4}>0,$$
so that (\ref{th:C-below}) follows by integration.

\subparagraph{\emph{\textmd{Proof of estimate (\ref{th:Cep-below})}}}

As before, we begin by showing that
\begin{equation}
	|A^{1/2}\uep(t)|^{2}e^{2\mu\Cep(t)}\geq k_{5}>0
	\quad\quad
	\forall t\geq 0.
	\label{est:Auep-Cep}
\end{equation}

To this end, let $v_{0}$ be the projection of $u_{0}$ in the
eigenspace of $\mu$, let $\uepm(t)$ be defined as in
(\ref{defn:uepn-uepm}), and let $\yep(t):=[\uepm(t)]^{2}$.  
From (\ref{eqn:uepn-uepm}) we have that
$$\yep'(t)=-2\mu\cep(t)\yep(t)-
2\ep\,\uepm(t)\cdot\uepm''(t).$$

\penalty -3000

Since $\yep(0)=|v_{0}|^{2}$, integrating this differential equation 
we obtain that
\begin{equation}
	e^{2\mu\Cep(t)}\yep(t)=|v_{0}|^{2}-2\ep\int_{0}^{t}
	\uepm(s)\cdot\uepm''(s)\cdot e^{2\mu\Cep(s)}\,ds.
	\label{est:yep}
\end{equation}

In order to estimate the last term, we integrate by parts, and we 
find that
\begin{eqnarray*}
	\lefteqn{\hspace{-3em}
	\int_{0}^{t} \uepm(s)\cdot\uepm''(s)\cdot
	e^{2\mu\Cep(s)}\,ds \ =\  \uepm(t)\cdot\uepm'(t)\cdot
	e^{2\mu\Cep(t)}- \uepm(0)\cdot\uepm'(0)}
	  \\
	  \noalign{\vspace{1ex}}
	  \hspace{3em} & & 
	  -\int_{0}^{t}[\uepm'(s)]^{2}e^{2\mu\Cep(s)}\,ds -\int_{0}^{t}
	 \uepm(s)\cdot\uepm'(s)\cdot 2\mu\cep(s)e^{2\mu\Cep(s)}\,ds
	 \\
	 \noalign{\vspace{1ex}}
	 & =: & A_{1}+A_{2}+A_{3}+A_{4}.
\end{eqnarray*}

Let us estimate the four terms separately. From (\ref{th:uepm}) we 
have that
\begin{equation}
	|\uepm(t)|\leq k_{6}e^{-\mu\Cep(t)}
	\quad\quad
	\forall t\geq 0,
	\label{est:uepm}
\end{equation}
\begin{equation}
	|\uepm'(t)|\leq k_{7}\cep(t)e^{-\mu\Cep(t)}
	\quad\quad
	\forall t\geq 0.
	\label{est:uepm'}
\end{equation}

Therefore, the estimate from above in (\ref{th:cep}) implies that
$$A_{1}\leq|\uepm(t)|\cdot|\uepm'(t)|\cdot e^{2\mu\Cep(t)}\leq 
k_{8}.$$

The term $A_{2}$ is a constant, independent of $\ep$, and of course 
$A_{3}\leq 0$.  

As for $A_{4}$, exploiting once more (\ref{est:uepm}) and
(\ref{est:uepm'}), we have that
$$-\uepm(t)\cdot\uepm'(t)\leq
|\uepm(t)|\cdot|\uepm'(t)|\leq
k_{9}\cep(t)e^{-2\mu\Cep(t)}.$$

Therefore, the estimate from above in (\ref{th:cep}) implies that
$$A_{4}\leq k_{10}\int_{0}^{t}[\cep(s)]^{2}ds \leq
k_{11}\int_{0}^{t}\frac{1}{(1+s)^{2}}\,ds\leq k_{11}.$$

In conclusion, we have proved that $A_{1}+A_{2}+A_{3}+A_{4}\leq 
k_{12}$. Coming back to (\ref{est:yep}), we have obtained that
$$e^{2\mu\Cep(t)}\yep(t)\geq|v_{0}|^{2}-2k_{12}\ep
\geq\frac{|v_{0}|^{2}}{2}
\quad\quad
\forall t\geq 0,$$
provided that $\ep$ is small enough. Therefore we have that
$$|A^{1/2}\uep(t)|^{2}e^{2\mu\Cep(t)}\geq 
\mu|\uepm(t)|^{2}e^{2\mu\Cep(t)}\geq\mu\frac{|v_{0}|^{2}}{2},$$
which implies (\ref{est:Auep-Cep}). It follows that
$$\left[e^{2\mu\gamma\Cep(t)}\right]'=
2\mu\gamma\cep(t)e^{2\mu\gamma\Cep(t)}=
2\mu\gamma\left(|A^{1/2}\uep(t)|^{2}e^{2\mu\Cep(t)}\right)^{\gamma}
\geq k_{13}>0$$
so that (\ref{th:Cep-below}) follows by integration.
\qed
\medskip

Next result is an estimate for a supersolution of an ordinary 
differential equation.

\begin{lemma}\label{lemma:ODE}
	Let $0<\delta\leq 1$ and $K>0$ be two real numbers.
	
	Then there exist $\ep_{1}>0$ and $M_{28}>0$, both depending on $\delta$
	and $K$, such that the following property holds true.  For every
	$\ep\in(0,\ep_{1})$, and every function $\psiep\in
	C^{1}([0,+\infty),\re)$ such that $\psiep(0)=1$, and
	\begin{equation}
		\psiep'(t)\geq -K\psiep(t)\left(
		\ep^{2\gamma}[\psiep(t)]^{2\gamma}+
		[\psiep(t)]^{2\gamma/\delta}\right)
		\quad\quad
		\forall t\geq 0,
		\label{hp:gep}
	\end{equation}
	we have that
	\begin{equation}
		\psiep\left(\frac{1}{\ep^{\delta}}\right)\geq 
		M_{28}\ep^{\delta^{2}/(2\gamma)}.
		\label{th:gep}
	\end{equation}
\end{lemma}

\paragraph{\emph{\textmd{Proof}}}

Let us consider the differential equation
\begin{equation}
	y'=-Ky\left\{\ep^{2\gamma}y^{2\gamma}+y^{2\gamma/\delta}\right\}.
	\label{ODE}
\end{equation}

Assumption (\ref{hp:gep}) is equivalent to saying that $\psiep(t)$ is a 
supersolution of (\ref{ODE}) for every $t\geq 0$. Let us set
$$z(t):=\left(\frac{\delta}{4K\gamma 
t+\delta}\right)^{\delta/(2\gamma)}
\quad\quad
\forall t\geq 0.$$

We claim that, for $\ep$ small enough, $z(t)$ is a subsolution of 
(\ref{ODE}) for $t\in[0,1/\ep^{\delta}]$. Indeed this is equivalent to 
showing that
$$z'(t)=-2K[z(t)]^{1+2\gamma/\delta}\leq
-Kz(t)\left\{ \ep^{2\gamma}[z(t)]^{2\gamma}+
[z(t)]^{2\gamma/\delta}\right\} 
\quad\quad 
\forall t\in[0,1/\ep^{\delta}],$$
which in turn is equivalent to
\begin{equation}
	[z(t)]^{(1-\delta)/\delta}\geq\ep
	\quad\quad 
	\forall t\in[0,1/\ep^{\delta}].
	\label{est:z}
\end{equation}

Since $z(t)$ is decreasing, and $0<\delta\leq 1$, it is enough to 
check (\ref{est:z}) when $t=1/\ep^{\delta}$. Now for $\ep$ small 
enough we have that
\begin{equation}
	\left[z\left(\frac{1}{\ep^{\delta}}\right)\right]^{(1-\delta)/\delta}
	=\left[\frac{\delta\ep^{\delta}}{4K\gamma+
	\delta\ep^{\delta}}\right]^{(1-\delta)/(2\gamma)}
	\geq k_{1}\ep^{\delta(1-\delta)/(2\gamma)}.
	\label{eqn:z}
\end{equation}

Since $\delta(1-\delta)/(2\gamma)< 1$, inequality (\ref{eqn:z})
implies (\ref{est:z}) when $\ep$ is small enough.

This proves that $z(t)$ is a subsolution of (\ref{ODE}) in the given 
interval. Since $z(0)=1=\psiep(0)$, the usual comparison principle 
yields that
$$\psiep\left(\frac{1}{\ep^{\delta}}\right)\geq
z\left(\frac{1}{\ep^{\delta}}\right)\geq
k_{2}\ep^{\delta^{2}/(2\gamma)},$$
which proves (\ref{th:gep}).
\qed

\subsection{Proof of Theorem~\ref{thm:main}}

This proof is organized as the proof of the main result
of~\cite{gg:de-dg1}.  The first part is the nonlinear core of the
argument, where we prove that $\roep(t)$ satisfies (\ref{th:roep}), 
of course under the assumption that $(u_{0},u_{1})\in 
D(A^{3/2})\times D(A^{1/2})$.
In the second part we apply Proposition~\ref{prop:linear} in order to
deduce all the conclusions of Theorem~\ref{thm:main}.

\subsubsection{Nonlinear core}

Let us set $\yep(t):=|\roep(t)|^{2}$. From (\ref{eqn:roep-lin}) and 
(\ref{defn:gep}) we have that
$$\yep'(t)=-2\langle\cep(t)A\uep(t)-c(t)Au(t),\roep(t)\rangle
-2\ep\langle\uep''(t),\roep(t)\rangle,$$
hence by (\ref{th:monoton}) it follows that
\begin{equation}
	\yep'(t)\leq-(\cep(t)+c(t))|A^{1/2}\roep(t)|^{2}
	-2\ep\langle\uep''(t),\roep(t)\rangle.
	\label{est:yep'1}
\end{equation}

Since $u(t)$ and $\uep(t)$, hence also $\roep(t)$, lie in the subspace
$H_{[\nu,+\infty)}$, we have that
$|A^{1/2}\roep(t)|^{2}\geq\nu|\roep(t)|^{2}$.  Thus (\ref{est:yep'1})
implies that 
$$\yep'(t)\leq-\nu(\cep(t)+c(t))\yep(t)
-2\ep\langle\uep''(t),\roep(t)\rangle.$$

Since $\yep(0)=0$, integrating this differential inequality we obtain
that
\begin{equation}
	\yep(t)\leq-2\ep e^{-\nu(\Cep(t)+C(t))}\int_{0}^{t}
	\langle\uep''(s),\roep(s)\rangle 
	e^{\nu(\Cep(s)+C(s))}\,ds.
	\label{est:yep-main}
\end{equation}

In order to estimate the last term, we integrate by parts, and we find 
that
\begin{eqnarray*}
	\lefteqn{\hspace{-6em}
	-2\ep\int_{0}^{t} \langle\uep''(s),\roep(s)\rangle
	e^{\nu(\Cep(s)+C(s))}\,ds \ =\  -2\ep
	\langle\uep'(t),\roep(t)\rangle e^{\nu(\Cep(t)+C(t))}} \\
	\hspace{3em} & & +2\ep \int_{0}^{t}\langle\uep'(s),\roep'(s)\rangle
	 e^{\nu(\Cep(s)+C(s))}\,ds \\
	 &  & +2\ep \int_{0}^{t} \langle\uep'(s),\roep(s)\rangle
	\nu(\cep(s)+c(s)) e^{\nu(\Cep(s)+C(s))}\,ds  \\
	\noalign{\vspace{1ex}}
	 & =: & A_{1}+A_{2}+A_{3}.
\end{eqnarray*}

Let us estimate the three terms separately.  From (\ref{th:uep'}) we
have that
$$2\ep |\uep'(t)|\cdot|\roep(t)|\leq
2\ep^{2}|\uep'(t)|^{2}+\frac{1}{2}|\roep(t)|^{2}\leq
\frac{2k_{1}\ep^{2}}{(1+t)^{2+1/\gamma}}
+\frac{1}{2}|\roep(t)|^{2},$$
hence
\begin{equation}
	A_{1}\leq 2\ep |\uep'(t)|\cdot|\roep(t)|\cdot
	e^{\nu(\Cep(t)+C(t))}\leq 
	\left(\frac{2k_{1}\ep^{2}}{(1+t)^{2+1/\gamma}}
	+\frac{1}{2}|\roep(t)|^{2}\right)
	e^{\nu(\Cep(t)+C(t))}.
	\label{est:A1}
\end{equation}

In order to estimate $A_{2}$, we first observe that
$$2\ep\langle\uep'(t),\roep'(t)\rangle\leq
2\ep|\uep'(t)|\cdot\left(|\rep'(t)|+|\tetep'(t)|\right)\leq
\ep^{2}|\uep'(t)|^{2}
+|\rep'(t)|^{2}+2\ep|\uep'(t)|\cdot|\tetep'(t)|.$$

Therefore, from (\ref{th:uep'}) and the explicit formula
(\ref{eqn:tetep}) for $\tetep(t)$, we obtain that
$$2\ep\langle\uep'(t),\roep'(t)\rangle\leq
k_{1}\frac{\ep^{2}}{(1+t)^{2+1/\gamma}}+|\rep'(t)|^{2}+
k_{2}\frac{\ep}{(1+t)^{1+1/(2\gamma)}}\cdot
e^{-t/\ep}.$$

On the other hand, from (\ref{th:C-above}), (\ref{th:Cep-above}), and
the fact that $\nu\leq\mu$, we have that
\begin{equation}
	e^{\nu(\Cep(t)+C(t))}\leq e^{\nu\Cep(t)}\cdot e^{\mu C(t)}\leq 
	k_{3}(1+t)^{1/\gamma},
	\label{est:exp-C-Cep}
\end{equation}
and therefore
\begin{eqnarray*}
	2\ep\langle\uep'(t),\roep'(t)\rangle e^{\nu(\Cep(t)+C(t))} & \leq & 
	k_{4}\frac{\ep^{2}}{(1+t)^{2}}+k_{3}|\rep'(t)|^{2}(1+t)^{1/\gamma}+
	k_{5}\frac{\ep}{(1+t)^{1-1/(2\gamma)}}e^{-t/\ep}  \\
	 & \leq & k_{4}\frac{\ep^{2}}{(1+t)^{2}}+k_{3}|\rep'(t)|^{2}(1+t)+
	 k_{5}\ep e^{-t/\ep}.
\end{eqnarray*}

Integrating in $[0,t]$, and exploiting (\ref{th:ee-rep'}), we obtain that
\begin{equation}
	A_{2}\leq k_{6}\ep^{2}.
	\label{est:A2}
\end{equation}

In order to estimate $A_{3}$, we first apply (\ref{th:c}),
(\ref{th:cep}), and (\ref{est:exp-C-Cep}) in order to obtain that
$$2\ep\langle\uep'(t),\roep(t)\rangle
\nu(\cep(t)+c(t)) e^{\nu(\Cep(t)+C(t))}\leq
k_{7}\ep|\uep'(t)|\cdot|\roep(t)|\cdot
\frac{1}{1+t}\cdot(1+t)^{1/\gamma}.$$

Now we estimate $|\uep'(t)|$ and $|\roep(t)|$ by means of
(\ref{th:uep'}) and (\ref{th:ee-roep}), respectively.  We obtain that
$$2\ep\langle\uep'(t),\roep(t)\rangle
\nu(\cep(t)+c(t)) e^{\nu(\Cep(t)+C(t))}\leq
k_{8}\frac{\ep^{2}}{(1+t)^{2-1/(2\gamma)}}.$$

Since $2-1/(2\gamma)>1$, integrating in $[0,t]$ we obtain that
\begin{equation}
	A_{3}\leq k_{9}\ep^{2}.
	\label{est:A3}
\end{equation}

Plugging (\ref{est:A1}), (\ref{est:A2}), and (\ref{est:A3}) into 
(\ref{est:yep-main}), we have that
\begin{equation}
	\yep(t)\leq k_{10}\frac{\ep^{2}}{(1+t)^{2+1/\gamma}}+
	\frac{1}{2}\yep(t)+k_{11}\ep^{2}e^{-\nu(\Cep(t)+C(t))}
	\label{est:yep1}
\end{equation}

Finally, (\ref{th:C-below}) and (\ref{th:Cep-below}) imply that
\begin{equation}
	e^{-\nu(\Cep(t)+C(t))}= 
	\left(e^{-\mu\Cep(t)}\cdot e^{-\mu C(t)}\right)^{\nu/\mu}\leq 
	k_{12}\frac{1}{(1+t)^{\delta/\gamma}}.
	\label{est:exp}
\end{equation}

Since $2+1/\gamma\geq 1/\gamma\geq \delta/\gamma$, plugging 
(\ref{est:exp}) into (\ref{est:yep1}) we find that
$$\frac{1}{2}\yep(t)\leq
k_{13}\frac{\ep^{2}}{(1+t)^{\delta/\gamma}},$$
which is equivalent to (\ref{th:roep}).

\subsubsection{Linear conclusion}

It remains to show that the assumptions of
Proposition~\ref{prop:linear} are satisfied. The argument is the same 
as in section~3.5 of~\cite{gg:de-dg1}, with the obvious changes in 
the exponents.

The a priori estimate on $|\roep(t)|$ is exactly the content of the
nonlinear core. The assumptions on $\cep(t)$ are exactly (\ref{th:cep}) and 
(\ref{th:cep'}).

In order to prove estimates on $\gep(t)$, we first estimate
$\cep(t)-c(t)$.  To this end, we apply the mean value theorem to the
function $\sigma^{\gamma}$, and we obtain the inequality
$$|y^{\gamma}-x^{\gamma}|\leq\gamma
\max\{y^{\gamma-1},x^{\gamma-1}\}\cdot|y-x| \quad\quad
\forall x\geq 0,\ \forall y\geq 0.$$

Setting $y:=|A^{1/2}\uep(t)|^{2}$ and 
$x:=|A^{1/2}u(t)|^{2}$, it follows that
\begin{equation}
	|\cep(t)-c(t)|\leq\gamma
	\max\left\{|A^{1/2}\uep|^{2(\gamma-1)},
	|A^{1/2}u|^{2(\gamma-1)}\right\}\cdot
	\left||A^{1/2}\uep|^{2}-|A^{1/2}u|^{2}\right|.
	\label{est:cep-c-1}
\end{equation}

From (\ref{th:Aqu}) and (\ref{th:Aquep}) we have that
\begin{equation}
	\max\left\{|A^{1/2}\uep(t)|^{2(\gamma-1)},
	|A^{1/2}u(t)|^{2(\gamma-1)}\right\}\leq
	\frac{k_{14}}{(1+t)^{1-1/\gamma}}.
	\label{est:cep-c-2}
\end{equation}

Moreover we have that
\begin{eqnarray*}
	\left||A^{1/2}\uep(t)|^{2}-|A^{1/2}u(t)|^{2}\right| & = & 
	\left|\langle A(\uep(t)+u(t)),\uep(t)-u(t)\rangle\right|\\
	 & \leq & \left(|A\uep(t)|+|Au(t)|\strut\right)
	 \cdot|\roep(t)|,
\end{eqnarray*}
so that from (\ref{th:Au}) and (\ref{th:Auep}) we obtain that
\begin{equation}
	\left||A^{1/2}\uep(t)|^{2}-|A^{1/2}u(t)|^{2}\right|\leq
	k_{15}\frac{|\roep(t)|}{(1+t)^{1/(2\gamma)}}.
	\label{est:cep-c-3}
\end{equation}

From (\ref{est:cep-c-1}), (\ref{est:cep-c-2}), (\ref{est:cep-c-3}), 
and (\ref{th:roep}), we conclude that
$$|\cep(t)-c(t)|\leq k_{16}
\frac{\ep}{(1+t)^{1+(\delta-1)/(2\gamma)}}.$$

From (\ref{th:Au}) we have therefore that
\begin{eqnarray}
	|\gep(t)|^{2} & \leq & 2(\cep(t)-c(t))^{2}|Au(t)|^{2}+
	2\ep^{2}|u''(t)|^{2}
	\nonumber  \\
	 & \leq & 
	 k_{17}\frac{\ep^{2}}{(1+t)^{2+\delta/\gamma}}+
	 2\ep^{2}|u''(t)|^{2}.
	\label{est:gep-main}
\end{eqnarray}

At this point (\ref{hp:gep-1}) follows from (\ref{th:u''-int}).
Moreover (\ref{hp:gep-2}) follows in an analogous way exploiting
(\ref{th:Atu}) instead of (\ref{th:Au}), and (\ref{th:auq-u''-int})
instead of (\ref{th:u''-int}).

Finally, from (\ref{est:gep-main}) and (\ref{th:u''}) we obtain 
that
$$|\gep(t)|^{2}\leq k_{17}
\frac{\ep^{2}}{(1+t)^{2+\delta/\gamma}}+
k_{18}\frac{\ep^{2}}{(1+t)^{4+1/\gamma}}\leq
k_{19}\frac{\ep^{2}}{(1+t)^{2+\delta/\gamma}},$$
where in the last inequality we used that $2+\delta/\gamma\leq 
4+1/\gamma$. This proves~(\ref{hp:gep-3}), and completes the proof of 
Theorem~\ref{thm:main}.\qed

\subsection{Proof of Theorem~\ref{thm:optimality}}

The assumptions on initial data guarantee that the solution $\uep(t)$
of the hyperbolic problem has only two Fourier components, and can be
written in the form
$$\uep(t):=\uepn(t)\frac{u_{1}}{|u_{1}|}+
\uepm(t)\frac{u_{0}}{|u_{0}|},$$
where the coefficients $\uepn(t)$ and $\uepm(t)$ are given by 
(\ref{defn:uepn-uepm}) (note that in this case $u_{0}=v_{0}$ and 
$u_{1}=v_{1}$).

On the other hand, the solution of the parabolic problem has only the
Fourier component with respect to $u_{0}$.  Since $u_{0}$ and $u_{1}$
are orthogonal, we can estimate the norm of $\roep(t)$ with the 
absolute value of its component with respect to $u_{1}$, namely
$$|\roep(t)|^{2}\geq\frac{\langle\roep(t),u_{1}\rangle^{2}}{|u_{1}|^{2}}=
\frac{\langle\uep(t),u_{1}\rangle^{2}}{|u_{1}|^{2}}
=|\uepn(t)|^{2}.$$

Therefore (\ref{th:optimality}) is proved if we show that
\begin{equation}
	\sup_{t\geq 0}\left\{(1+t)^{\delta/\gamma}|\uepn(t)|^{2}\right\}\geq 
	k_{1}\ep^{2}.
	\label{th:optimality-red}
\end{equation}

\paragraph{\textmd{\emph{Estimate of $\uepn$ from above}}}

We prove that
\begin{equation}
	|\uepn(t)|\leq k_{2}\ep e^{-\nu\Cep(t)}
	\quad\quad
	\forall t\geq 0.
	\label{th:uepn-above}
\end{equation}

Let us consider in (\ref{eqn:uepn-uepm}) the equation solved by
$\uepn(t)$.  Moving $\ep\uepn''(t)$ to the right-hand side, we can
interpret it as a first order equation.  Since $\uepn(0)=0$,
integrating this differential equation we obtain that
$$\uepn(t)=-\ep e^{-\nu\Cep(t)}\int_{0}^{t}
\uepn''(s)e^{\nu\Cep(s)}\,ds.$$

Integrating by parts, and remarking that $\uepn'(0)=|u_{1}|$, this 
can be rewritten as
\begin{equation}
	\uepn(t)=-\ep\uepn'(t)+\ep |u_{1}|e^{-\nu\Cep(t)}+
	\ep e^{-\nu\Cep(t)}\int_{0}^{t}
	\uepn'(s)\nu\cep(s)e^{\nu\Cep(s)}\,ds,
	\label{est:uepn-parts}
\end{equation}
hence
$$|\uepn(t)|\leq\ep e^{-\nu\Cep(t)}\left\{
|\uepn'(t)|\cdot e^{\nu\Cep(t)}+|u_{1}|+\nu\int_{0}^{t} 
|\uepn'(s)|\cdot\cep(s)e^{\nu\Cep(s)}\,ds\right\}.$$

Now from (\ref{th:uepn}), and the estimate from above in (\ref{th:cep}), 
we have that
\begin{equation}
	|\uepn'(t)|\cdot e^{\nu\Cep(t)}\leq
	k_{3}\cep(t)\leq\frac{k_{4}}{1+t},
	\label{est:uepn'}
\end{equation}
hence
$$|\uepn(t)|\leq\ep e^{-\nu\Cep(t)}\left\{
\frac{k_{4}}{1+t}+|u_{1}|+k_{5}\int_{0}^{t} 
\frac{1}{(1+s)^{2}}\,ds\right\},$$
which easily implies (\ref{th:uepn-above}).

\paragraph{\textmd{\emph{Estimate of $\uepn$ from below}}}

We prove that
\begin{equation}
	\uepn(t)\geq \ep e^{-\nu\Cep(t)}\left\{
	|u_{1}|-\frac{k_{6}}{1+t}-k_{7}\ep\right\}
	\quad\quad
	\forall t\geq 0.
	\label{th:uepn-below}
\end{equation}

To this end, we need to estimate the absolute value of the last term
in (\ref{est:uepn-parts}).  Thus we first integrate by parts, and we
obtain that
\begin{eqnarray*}
	\int_{0}^{t} \uepn'(s)\nu\cep(s)e^{\nu\Cep(s)}\,ds & = &
	\uepn(t)\nu\cep(t)e^{\nu\Cep(t)}
	-\int_{0}^{t} \uepn(s)\nu\cep'(s)e^{\nu\Cep(s)}\,ds \\
	 &  & 
	-\int_{0}^{t} \uepn(s)\nu^{2}\cep^{2}(s)e^{\nu\Cep(s)}\,ds \\
	\noalign{\vspace{1ex}}
	 & =: & A_{1}+A_{2}+A_{3}.
\end{eqnarray*}

From (\ref{th:uepn-above}), and the estimate from above in 
(\ref{th:cep}), we have that
\begin{equation}
	|A_{1}|\leq k_{8}\ep.
	\label{est:A1-opt}
\end{equation}

In order to control $A_{2}$ and $A_{3}$, we estimate $|\uepn(s)|$ by 
means of (\ref{th:uepn-above}), we estimate $\cep(s)$ by means of 
(\ref{th:cep}), and $|\cep'(s)|$ by exploiting (\ref{th:cep}) and 
(\ref{th:cep'}) in order to deduce that
$$|\cep'(s)|=\frac{|\cep'(s)|}{\cep(s)}\cdot\cep(s)\leq
\frac{k_{9}}{(1+s)^{2}}.$$

Thus we obtain that
\begin{equation}
	|A_{2}|+|A_{3}|\leq k_{10}\ep\int_{0}^{t}\frac{1}{(1+s)^{2}}\,ds
	\leq k_{10}\ep.
	\label{est:A2+A3}
\end{equation}

Plugging (\ref{est:A1-opt}) and (\ref{est:A2+A3}) into 
(\ref{est:uepn-parts}), and exploiting (\ref{est:uepn'}), we obtain that
\begin{eqnarray*}
	\uepn(t) & \geq & -\ep|\uepn'(t)|+\ep|u_{1}| e^{-\nu\Cep(t)}
	-k_{11}\ep^{2}e^{-\nu\Cep(t)}
	  \\
	 \noalign{\vspace{1ex}} 
	 & = & \ep e^{-\nu\Cep(t)}\left\{|u_{1}|-
	 |\uepn'(t)|e^{\nu\Cep(t)}-k_{11}\ep\right\}   \\
	 \noalign{\vspace{1ex}} 
	  & \geq & \ep e^{-\nu\Cep(t)}\left\{
	|u_{1}|-\frac{k_{4}}{1+t}-k_{11}\ep\right\},
\end{eqnarray*}
which is exactly (\ref{th:uepn-below}).

\paragraph{\textmd{\emph{Estimate of the exponential}}}

Let us set for simplicity $\psiep(t):=e^{-\nu\Cep(t)}$. We prove that
\begin{equation}
	\psiep\left(\frac{1}{\ep^{\delta}}\right)\geq 
	k_{12}\ep^{\delta^{2}/(2\gamma)}.
	\label{est:gep}
\end{equation}

Indeed from (\ref{th:uepn-above}) and (\ref{th:uepm}) we have that
$$|A^{1/2}\uep(t)|^{2}=\nu|\uepn(t)|^{2}+\mu|\uepm(t)|^{2}\leq
k_{13}\left(\ep^{2}[\psiep(t)]^{2}+[\psiep(t)]^{2\mu/\nu}\right),$$
hence
$$\cep(t)=|A^{1/2}\uep(t)|^{2\gamma}\leq k_{14}
\left(\ep^{2\gamma}[\psiep(t)]^{2\gamma}
+[\psiep(t)]^{2\gamma/\delta}\right).$$

Since $\psiep'(t)=-\nu\cep(t)\psiep(t)$, this implies that
$$\psiep'(t)\geq -k_{15}\psiep(t)
\left(\ep^{2\gamma}[\psiep(t)]^{2\gamma}
+[\psiep(t)]^{2\gamma/\delta}\right)
\quad\quad
\forall t\geq 0.$$

At this point (\ref{est:gep}) follows from Lemma~\ref{lemma:ODE}.

\paragraph{\textmd{\emph{Conclusion}}}

We are now ready to prove (\ref{th:optimality-red}). Let us set 
$t:=1/\ep^{\delta}$ in (\ref{th:uepn-below}). We obtain that
$$\uepn\left(\frac{1}{\ep^{\delta}}\right)\geq
\ep\psiep\left(\frac{1}{\ep^{\delta}}\right)
\left\{|u_{1}|-k_{6}\frac{\ep^{\delta}}{1+\ep^{\delta}}
-k_{7}\ep\right\}.$$

If $\ep$ is small enough, the term between braces is larger than or 
equal to $|u_{1}|/2$, hence by (\ref{est:gep}) we obtain that
$$\uepn\left(\frac{1}{\ep^{\delta}}\right)\geq
\frac{|u_{1}|}{2}\ep\psiep\left(\frac{1}{\ep^{\delta}}\right)\geq
k_{16}\ep\cdot\ep^{\delta^{2}/(2\gamma)}.$$

Therefore we conclude that
$$\sup_{t\geq 0}|\uepn(t)|^{2}(1+t)^{\delta/\gamma}\geq
\left|\uepn\left(\frac{1}{\ep^{\delta}}\right)\right|^{2}
\left(\frac{1+\ep^{\delta}}{\ep^{\delta}}\right)^{\delta/\gamma}\geq
k_{17}\ep^{2},$$
which completes the proof.
\qed

\label{NumeroPagine}


\begin{thebibliography}{99}

	\bibitem{ch}{\sc R.\ Chill, A.\ Haraux}; An optimal estimate for the
	time singular limit of an abstract wave equation.  \emph{Funkcial.\  Ekvac.}\
	\textbf{47} (2004),  no.~2, 277--290.

	\bibitem{ghisi:error}{\sc M.\ Ghisi}; Hyperbolic-parabolic
	singular perturbation for mildly degenerate Kirchhoff equations
	with weak dissipation.  \emph{Adv.\ Differential Equations}
	\textbf{17} (2012), no.~1--2, 1--36.

	\bibitem{ghisi:decay}{\sc M.\ Ghisi}; Asymptotic limits for mildly
	degenerate Kirchhoff equations.  Preprint. {\tt arXiv:1105.5358v1
	[math.AP]}.

	\bibitem{gg:l-cattaneo}{\sc M.\ Ghisi, M.\ Gobbino}; Global-in-time
	uniform convergence for linear hyperbolic-parabolic singular
	perturbations.  \emph{Acta Math.\ Sin.\ (Engl.\  Ser.)} \textbf{22}
	(2006), no.~4, 1161--1170.

	\bibitem{gg:k-decay}{\sc M.\ Ghisi, M.\ Gobbino}; Hyperbolic-parabolic
	singular perturbation for mildly degenerate Kirchhoff equations:
	time-decay estimates.  \emph{J.\ Differential Equations} \textbf{245}
	(2008), no.~10, 2979--3007.

	\bibitem{gg:k-PS}{\sc M.\ Ghisi, M.\ Gobbino};
	Hyperbolic-parabolic singular perturbation for mildly degenerate
	Kirchhoff equations: global-in-time error estimates.
	\emph{Commun.\ Pure Appl.\ Anal.}\ \textbf{8} (2009), no.~4,
	1313--1332.

	\bibitem{gg:w-ndg}{\sc M.\ Ghisi, M.\ Gobbino};
	Hyperbolic-parabolic singular perturbation for nondegenerate
	Kirchhoff equations with critical weak dissipation.  \emph{Math.\
	Ann.}\ doi: 10.1007/s00208-011-0765-x.

	\bibitem{gg:survey-diss}{\sc M.\ Ghisi, M.\ Gobbino};
	Hyperbolic-parabolic singular perturbation for Kirchhoff equations
	with weak dissipation.  \emph{Rend.\ Ist.\ Mat.\ Univ.\ Trieste}
	\textbf{42} Suppl.\ (2010), 67--88.
	
	\bibitem{gg:de-dg1}{\sc M.\ Ghisi, M.\ Gobbino};
	Hyperbolic-parabolic singular perturbation for mildly degenerate
	Kirchhoff equations: Decay-error estimates.  \emph{J.\ Differential
	Equations} (2012), doi:10.1016/j.jde.2012.02.019.

	\bibitem{k-par}{\sc M.\ Gobbino}; Quasilinear degenerate parabolic
	equations of Kirchhoff type.  \emph{Math.\  Methods Appl.\  Sci.}\ 
	\textbf{22} (1999), no.~5, 375--388.

	\bibitem{haraux}\textsc{A.\ Haraux}; Slow and fast decay of
	solutions to some second order evolution equations.  \emph{J.\
	Anal.\ Math.}\ \textbf{95} (2005), 297--321.
	
	\bibitem{yamazaki} {\sc H.\ Hashimoto, T.\ Yamazaki};
	Hyperbolic-parabolic singular perturbation for quasilinear equations
	of Kirchhoff type.  \emph{J.\ Differential Equations} \textbf{237}
	(2007), no.~2, 491--525.

	\bibitem{lions}{\sc J. L. Lions}; {\em Perturbations singuli\'{e}res
	dans les probl\`{e}mes aux limites et en control optimal}, Lecture
	Notes in Mathematics, Vol.\  323.  Springer-Verlag, Berlin-New York,
	1973.

	\bibitem{ny}{\sc K.\ Nishihara, Y.\ Yamada}; On global solutions of
	some degenerate quasilinear hyperbolic equations with dissipative
	terms.  \emph{Funkcial.\  Ekvac.}\ \textbf{33} (1990), no.~1, 151--159.

	\bibitem{rudin} {\sc W. Rudin}; {\em Functional Analysis}, 
	McGraw-Hill, New York, 1973.

	\bibitem{yamazaki-wd}{\sc T.\ Yamazaki}; Hyperbolic-parabolic
	singular perturbation for quasilinear equations of Kirchhoff type
	with weak dissipation.  \emph{Math.\ Methods Appl.\ Sci.}\
	\textbf{32} (2009), no.~15, 1893--1918.

	\bibitem{yamazaki-cwd}{\sc T.\ Yamazaki}; Hyperbolic-parabolic
	singular perturbation for quasilinear equations of Kirchhoff type with
	weak dissipation of critical power.  Preprint. 

\end{thebibliography}
\end{document}